\documentclass[11pt]{article}
\usepackage{microtype}

\usepackage{amsmath,amssymb}
\usepackage{amsfonts,bm,bbm} 
\usepackage{amssymb,amsthm,enumitem}
\usepackage{tikz}
\usetikzlibrary{decorations.pathmorphing,positioning,patterns}
\usepackage{multirow}
\usepackage{xcolor,graphicx,float,booktabs}
\usepackage{tabularx}
\newcolumntype{L}[1]{>{\raggedright\arraybackslash}m{#1}}
\newcolumntype{C}[1]{>{\centering\arraybackslash}m{#1}}
\newcolumntype{R}[1]{>{\raggedleft\arraybackslash}m{#1}}

\usepackage{algorithm}
\usepackage{algorithmic}

\usepackage{epstopdf}

\usepackage{verbatim} 
\allowdisplaybreaks[4]
\usepackage{authblk}

\usepackage{hyperref}
\hypersetup{colorlinks=true,
linkcolor=blue,
anchorcolor=blue,
citecolor=blue}

\usepackage[round]{natbib}
\usepackage[
margin=1in,
includefoot,
footskip=30pt,
]{geometry}

\newcommand{\cS}{\mathcal{S}}

\newcommand{\cL}{\mathcal{L}}
\newcommand{\cX}{\mathcal{X}}

\newcommand{\cI}{\mathcal{I}}
\newcommand{\cJ}{\mathcal{J}}
\newcommand{\bX}{\mathbf{X}}
\newcommand{\bx}{\mathbf{x}}
\newcommand{\E}{\mathbb{E}}
\newcommand{\R}{\mathbb{R}}
\renewcommand{\S}{\mathbb{S}}

\newcommand{\diam}{\mathrm{diam}}
\newcommand{\dist}{\mathrm{dist}}

\renewcommand{\P}{\mathbb{P}}
\newcommand{\1}{\mathbbm{1}}

\newtheorem{theorem}{Theorem}[section]
\newtheorem{lemma}[theorem]{Lemma}
\newtheorem{definition}[theorem]{Definition}

\newtheorem{assumption}{Assumption}[section]
\numberwithin{equation}{section}

\title{Regular Tree Search for Simulation Optimization}

\author[1, 2]{Du-Yi Wang\footnote{tylzml@ruc.edu.cn}}
\author[1]{Guo Liang\footnote{liangguo000221@ruc.edu.cn}}
\author[2]{Guangwu Liu\footnote{msgw.liu@cityu.edu.hk}}
\author[1]{Kun Zhang\footnote{kunzhang@ruc.edu.cn}}
\affil[1]{Institute of Statistics and Big Data\protect\\Renmin University of China\protect\\Beijing, China}
\affil[2]{Department of Decision Analytics and Operations\protect\\City University of Hong Kong\protect\\Tat Chee Avenue, Kowloon, Hong Kong, China}

\begin{document}
\normalsize

\maketitle

\begin{abstract}

Tackling simulation optimization problems with non-convex objective functions remains a fundamental challenge in operations research. In this paper, we propose a class of random search algorithms, called Regular Tree Search, which integrates adaptive sampling with recursive partitioning of the search space. The algorithm concentrates simulations on increasingly promising regions by iteratively refining a tree structure. A tree search strategy guides sampling decisions, while partitioning is triggered when the number of samples in a leaf node exceeds a threshold that depends on its depth. Furthermore, a specific tree search strategy, Upper Confidence Bounds applied to Trees (UCT), is employed in the Regular Tree Search. We prove global convergence under sub-Gaussian noise, based on assumptions involving the optimality gap, without requiring continuity of the objective function. Numerical experiments confirm that the algorithm reliably identifies the global optimum and provides accurate estimates of its objective value.

\end{abstract}

\section{Introduction}
\label{sec:intro}
Simulation optimization, also known as optimization via simulation (OvS), refers to a class of optimization problems where the objective function values can only be evaluated through simulation, and the evaluations are subject to stochastic noise. Such problems have wide applications in various domains, including engineering design and hyperparameter tuning. For a comprehensive review of simulation optimization, see \cite{amaran2016simulation,fan2024review}.

When the decision variables are discrete, the problem is called Discrete OvS (DOvS). Many algorithms for DOvS are based on random search strategies, such as the COMPASS algorithm \citep{hongDiscreteOptimizationSimulation2006}. When the decision space is finite, the problem can be formulated as a multi-armed bandit problem, with well-established methods like Upper Confidence Bound (UCB) aiming to minimize regret \citep{auer2002finitetime}. Alternatively, when the goal is to accurately estimate the optimal value rather than minimize regret, algorithms such as those proposed by \cite{9004921,kun2024upperconfidenceboundapproach} are more appropriate. If the goal is to identify the best alternative (or arm), the problem is generally referred to as a Ranking and Selection problem. For comprehensive reviews of this setting, see \cite{Chen2015,hong2021review}.

When the decision variables are continuous, the problem is known as Continuous OvS (COvS). 
A wide range of algorithms has been developed to address COvS problems. Among these are stochastic approximation algorithms \citep{kiefer1952stochastic,spall1992multivariate,doi:10.1142/S0217595925400056}, response surface methodologies \citep{chang2013stochastic}, model-based methods \citep{hu2007model}, and random search algorithms \citep{shi2000nested,andradottir2010adaptive,kiatsupaibulSingleObservationAdaptive2018,zhang2022actorcriticlikea,wang2025gaussian}.
Different algorithms are suited to different COvS problem settings and offer distinct convergence guarantees. Notably, many algorithms assume global continuity of the objective function, an assumption that may fail in practice.
In this paper, our focus is on random search algorithms, with the goal of designing methods that are gradient-free, globally convergent, and capable of handling discontinuous objective functions.

We consider a class of random search algorithms that partition the search space into disjoint regions to explicitly resolve the exploration-exploitation dilemma. This structural decomposition creates isolated promising areas, enabling strategic navigation of the fundamental tension between exploring uncertain regions and exploiting promising ones. The resulting partition structure enables discrete optimization algorithms to be seamlessly embedded within this framework, thereby transforming continuous optimization problems into the identification of optimal subspaces within discretized regions. Tree structures naturally embody this paradigm, as they recursively split the domain according to criteria that progressively refine promising regions.
One of the most well-known decision tree algorithms is Classification And Regression Tree (CART), introduced by \cite{breiman1984classification}. At each node, CART selects axis-aligned split directions to optimize the Gini index for classification or mean squared error (MSE) for regression. Subsequent developments, including random forests \citep{breimanRandomForests2001} and various theoretical advancements \citep{scornet2015consistency,10.5555/3666122.3668639,wagerEstimationInferenceHeterogeneous2018,zhang2024adaptive}, have improved predictive performance and provided convergence guarantees for these models. However, despite their contributions to predictive modeling, these advances often fail to address the unique challenges of simulation optimization, where the ultimate goal is to identify the global optimum.

A key distinction between prediction and simulation optimization lies in sample usage: prediction utilizes a complete sample set, whereas simulation optimization depends on samples generated iteratively. This fundamental distinction raises a critical question: which point should be simulated next? A central challenge in this setting is selecting the leaf node within a hierarchical tree structure, a key element in effective space partitioning. In addressing this challenge, Monte Carlo Tree Search (MCTS) provides a principled framework that adaptively balances exploration and exploitation. For instance, the UCT algorithm proposed by \cite{changAdaptiveSamplingAlgorithm2005,10.1007/11871842_29} exemplifies this approach. The practical success of UCT-based MCTS in complex combinatorial environments, as demonstrated by systems like AlphaGo \citep{silver2016mastering}, highlights the potential of tree-based strategies for tackling challenging optimization problems.

Specialized algorithms have been developed to adapt tree structures for optimization tasks. The optimistic optimization algorithm for noise-free settings proposed by \cite{munosOptimisticOptimizationDeterministic2011} and its extensions to noisy settings \citep{valko2013stochastic} rely on tree-based partitioning that employs random directional selection combined with uniform splitting. However, due to their reliance on uniform partitioning, these algorithms are not well-suited for adaptive splitting and tend to be less effective in isolating promising regions. 
In contrast, adaptive partitioning strategies have been introduced for noise-free settings. These include using linear regressors for space partitioning in neural architecture search \citep{wang2021sampleefficient}, and distinguishing between high- and low-value regions via sample-driven splits \citep{10.5555/3495724.3497361}. Despite promising experimental results, these methods lack theoretical guarantees.

The remainder of this paper is organized as follows. In Section \ref{sec:problem_formulation}, we formulate the simulation optimization problem and outline the key assumptions. Section \ref{sec:regular_tree} introduces the concept of the Regular Tree and details the criterion for its construction. In Section \ref{sec:regular_tree_search}, we describe the Regular Tree Search and explain how the UCT algorithm is integrated into it. Section \ref{sec:convergence_analysis} presents the theoretical foundations that guarantee the global convergence of the algorithms under certain conditions. Finally, Section \ref{sec:numerical_experiments} presents numerical results that validate the proposed algorithm, and Section \ref{sec:conclusion} concludes the paper.

\section{Problem Formulation}\label{sec:problem_formulation}
We consider the following simulation optimization problem:
\begin{align*}
    \max_{\bx\in \cX \subseteq \R^d} \{ \mu (\bx) = \E [Y| \bX=\bx] \}.
\end{align*}
Here, $\bx \in \R^d$ denotes the decision variable, $\cX$ denotes the feasible domain, $\mu(\bx)$ denotes the unknown response surface of the simulation model, and $Y$ represents the random output of a simulation model evaluated at $\bx$. Since the conditional distribution of $ Y | \bX = \bx$ is unknown and $\mu(\bx)$ has no closed form, we rely on running simulation experiments to generate independent samples of $Y|\bX = \bx$, denoted by $Y(\bx)$:
\begin{align*}
    Y(\bx) = \mu(\bx) + \epsilon(\bx),
\end{align*}
where $\epsilon(\bx)$ is the zero mean simulation noise. 

Our analysis focuses on the setting where the feasible domain $\cX$ is a compact subset of $\R^d$. While our theoretical framework can be adaptable to more general compact sets (as any compact set can be enclosed in a sufficiently large hyperrectangle), we specifically consider the canonical unit hypercube $\cX = [0, 1]^d$ for analytical tractability. 

A fundamental challenge in simulation optimization arises from the black box of $\mu$. Algorithms with minimal structural assumptions are highly desirable. This enhances the algorithm's applicability to a broad range of real-world problems where prior knowledge of the objective landscape is limited or unavailable.
Many established algorithms, particularly those within the frameworks of Bayesian optimization and methods ensuring global convergence, rely on the assumption of global continuity. Examples include the approaches discussed in \cite{frazier2018tutorial} and \cite{kiatsupaibulSingleObservationAdaptive2018}.
Although this assumption facilitates theoretical analysis, it is restrictive and may be violated in practical simulation settings. For instance, functions with discontinuities might arise in specific system models. 
We consider the conditions focused on the optimal point instead of global continuity.

Throughout this paper, we make the following assumptions:
\begin{assumption}
\label{assum:bound}
    The simulation noise $\epsilon(\bx)$ is sub-Gaussian with the variance proxy $C_{\text{sg}}^2$ for all $\bx \in \cX$ , i.e., $\mathbb{E}\left[e^{\lambda \epsilon(\bx)}\right] \leq e^{\lambda^2 C_{\text{sg}}^2/2}, \forall \lambda \in \mathbb{R}.$
\end{assumption}
\begin{assumption}
\label{assum:local_increase}
    There exists a unique optimal point $\bx^*\in \cX$ of $\mu$ satisfying $\mu(\bx^*) = \sup_{\bx\in \cX}\mu(\bx)$.
    The function $\mu$ exhibits local smoothness: there exists a constant $C>0$ such that for any $\bx\in \cX$:
    \begin{align*}
        \mu (\bx^*) - \mu (\bx) \leq C\cdot \| \bx- \bx^* \|.
    \end{align*}
    where $ \| \cdot \|$ denotes the Euclidean norm.
\end{assumption}
\begin{assumption}
\label{assum:local_decay}
    For any $\epsilon > 0$, there exists $\delta(\epsilon) > 0$ such that for all $\bx \in \cX$ satisfying $\| \bx- \bx^* \| \geq \epsilon$, 
    \begin{align*}
        \mu(\bx) < \mu(\bx^*) - \delta(\epsilon).
    \end{align*}
\end{assumption}

The sub-Gaussian noise assumption is a standard condition in the simulation optimization literature. Although the uniqueness condition in Assumption \ref{assum:local_increase} could be relaxed to accommodate multiple optimal points, we retain it to highlight our algorithmic innovation.

It is important to note that continuity does not imply Assumption \ref{assum:local_increase}. For example, the one-dimensional function $\mu(\bx) = 1-\sqrt{\bx}$ is continuous on $[0, 1]$ but has an infinite derivative at its optimal point $0$. Assumption \ref{assum:local_increase} is designed to capture the quantitative relationship between the optimal point and other points, without imposing a global smoothness condition.

Assumption \ref{assum:local_decay} is weaker than global continuity. To see this, suppose $\mu$ is continuous on $\cX$ with unique optimal point $\bx^*$. Then, for any $\epsilon>0$, the set $ K_\epsilon = \{ \bx \in \cX : \| \bx- \bx^* \| \geq \epsilon \}$ is compact by the boundedness and closeness. By the Extreme Value Theorem, the continuous function $\mu(\bx)$ attains its maximum on the compact set $K_\epsilon$, denoted by $M_\epsilon = \max_{\bx \in K_\epsilon} \mu(\bx)$.
Since $\bx^*$ is the unique global maximum of $\mu(\bx)$ and $\bx^* \notin K_\epsilon$, it follows that $M_\epsilon < \mu(\bx^*)$. Taking $\delta(\epsilon) = \frac{\mu(\bx^*)-M_\epsilon}{2}>0$, then for any $\bx \in K_\epsilon$, we have:
\begin{align*}
    \mu(\bx) \leq M_\epsilon = \mu(\bx^*) - 2\delta(\epsilon) < \mu(\bx^*) - \delta(\epsilon),
\end{align*}
which directly verifies Assumption \ref{assum:local_decay}. 
Hence, Assumption \ref{assum:local_decay} describes a local property of the optimum and is strictly weaker than global continuity. Assumption \ref{assum:local_increase} captures how the function rises toward the optimum, while Assumption \ref{assum:local_decay} ensures it declines beyond a neighborhood. Both express local characteristics of the optimum without relying on global smoothness.

\section{Regular Tree}\label{sec:regular_tree}
To address the theoretical limitations of the original CART algorithm, particularly the lack of convergence guarantees, \cite{wagerEstimationInferenceHeterogeneous2018} introduce a set of definitions that ensure favorable statistical properties. In this work, we adapt and extend these definitions into four definitions, honesty, random‑split, $\alpha$‑spatial balance, and $f(c)$‑sample balance, to align with the specific requirements of simulation optimization. 

\begin{definition}
\label{def:honest}
    A tree is \textbf{honest} if, for each training sample $i$, it only uses the response $Y_i$ to estimate the value of the leaf or to decide where to place the splits, but not both.
\end{definition}
The concept of honesty ensures that the data used to grow the tree is entirely independent of the data used to estimate the leaf values. Although this separation may lose efficiency since it prevents using each observation for both purposes, it guarantees that the sample mean estimated in each leaf node is an unbiased estimator of the true expected value. 
In practice, we implement honesty by maintaining two disjoint datasets, $\S_\cI$ and $\S_\cJ$. The dataset $\S_\cI$ is used to estimate the values at the leaf nodes, while $\S_\cJ$ is used to grow the tree (i.e., to determine the split direction and split value). 
Each leaf node $L$ corresponds to a hyperrectangle in the feasible domain.
We define the random variable $Y(L)$ through the following two-step sampling procedure:
\begin{enumerate}
    \item Sample $\bx$ uniformly from $L$.
    \item Generate $Y$ according to the conditional distribution $\P_{Y|\bX=\bx}$.
\end{enumerate}

Let $\bar{Y}_L = \frac{\sum_{i\in \cI} \1{\{ \bx_i\in L \}} Y_i}{ \sum_{i\in \cI} \1{ \{ \bx_i\in L \} } }$ denote the sample mean of responses in $L$ computed from $\S_\cI$. By construction, we have
\begin{align*}
    \E[\bar{Y}_L] = \E[Y (L)].
\end{align*}

\begin{definition}
\label{def:random_split}
    A tree is \textbf{random-split} if, at each step of the tree-growing procedure, the probability that the next split occurs along the $j$-th feature is bounded below by $\kappa/d$ for some $0< \kappa \leq 1$, for all $j=1, \dots, d$.
\end{definition} 

\begin{definition}
\label{def:spatially_balanced}
    A tree grown by recursive partitioning is \textbf{$\alpha$-spatially balanced} if, after each split, the Lebesgue measure of the splitting coordinate in each child node is at least an $\alpha$ fraction of that in the parent node.
\end{definition}
Random split ensures that every feature has a non-negligible chance of being selected at every split. By imposing a uniform lower bound $\kappa/d$ on the probability of selecting any given coordinate, we prevent the algorithm from systematically neglecting certain dimensions. This balance in feature selection is essential for deriving probabilistic bounds on how the tree explores the feasible domain and ensuring that all directions are sufficiently considered during the tree's growth.
$\alpha$‑spatial balance further controls the geometry of each split. This condition prevents splits from creating one tiny region and one huge region, instead guaranteeing that every split reduces the size of each subregion at a predictable rate.
Together, these two definitions ensure that the diameter of a leaf, defined as the Euclidean norm of its longest side, converges to 0 in probability as the leaf's depth increases. Intuitively, random‑split spreads the splits across all dimensions, while $\alpha$‑spatial balance prevents any single split from producing excessively imbalanced partitions. Specifically, we have the following result:
\begin{lemma}
\label{lemma:diam_depth}
    Let $T$ be a random-split, $\alpha$-spatially balanced tree, and let $L$ be a leaf of $T$ with depth $c$. Suppose that leaf $L = \bigotimes_{j=1}^d \left[ r_j^-, r_j^+ \right]$, where $0 \leq r_j^- < r_j^+ \leq 1$. Define the diameter of $L$ as $\diam(L) = \sqrt{\sum_{j=1}^d (r_j^+ - r_j^-)^2}$. 
    For any $\lambda\in (0, 1)$, the following inequality holds:
    \begin{align*}
        \P \left[\diam(L) \geq \sqrt{d} (1-\alpha)^{(1-\lambda)\frac{\kappa}{d}c}\right] \leq d \exp \left\{-\frac{\lambda^2}{2} \frac{\kappa}{d} c \right\}.
    \end{align*}
\end{lemma}

\begin{definition}
\label{def:sample_balanced}
    A tree is \textbf{$f(c)$-sample balanced} if the number of $\cI$-samples in the leaf node at depth $c$ is bounded between $\beta f(c-1)$ and $f(c)-1$ for a positive function $f(c)$ and $0 < \beta < 1/2$.
\end{definition}

The function $f(c)$ generalizes the concept of $k$-regularity from \cite{wager2015adaptive}, where $f(c)$ is always a fixed constant. 
In contrast to the original CART setting, where the entire dataset is provided in advance, simulation optimization involves data generated adaptively as the algorithm progresses. Due to the presence of noise in the observations, accurately identifying the optimal point and estimating the optimal value requires that the number of samples in each leaf increases as the algorithm progresses. Consequently, in our framework, the number of samples in a leaf is linked to the depth of the leaf: deeper leaf nodes are required to contain more samples to ensure statistical reliability.

A tree that satisfies all four definitions is called a \textit{Regular Tree}. To construct such a tree, we need to design a splitting criterion that explicitly enforces these four definitions. We propose a splitting criterion that guarantees compliance with all four definitions, as detailed in Algorithm \ref{algorithm:splitting_criterion}. 

\begin{algorithm}[htbp]
\caption{Regular Tree splitting criterion}
\renewcommand{\algorithmicrequire}{\textbf{Input:}}
\label{algorithm:splitting_criterion}
\begin{algorithmic}[1]
\REQUIRE Datasets $\S_\cI$ and $\S_\cJ$, spatially balanced parameter $\alpha \in (0, 0.5]$, random-split parameter $\kappa \in (0, 1]$, rectangle $L = \bigotimes_{j=1}^d \left[ r_j^-, r_j^+ \right]$ with depth $c$ where $0 \leq r_j^- < r_j^+ \leq 1$, and sample balanced parameter $f(c), \beta \in (0, 0.5)$ . 
\STATE \textbf{With probability $\kappa$:}
\STATE \quad Randomly select a direction $j^*\in \{1, 2, \dots, d\}$.
\STATE \quad Select the best split value $z^*$ by optimizing the MSE on $\S_\cJ$ along $j^*$-th direction:
\begin{align*}
    z^* = \arg\min_{z}\sum_{i \in \cJ} ( Y_i - \Bar{Y}_1 )^2 \1 \{ \bx_i \in L_1 \} + \sum_{i \in \cJ} (Y_i - \bar{Y}_2 )^2 \1 \{ \bx_i \in L_2 \},
\end{align*}
\quad ensuring:
\begin{align}
\label{splitting_criterion_1}
    (1-\alpha) r_{j^*}^- + \alpha r_{j^*}^+ \leq z^* \leq \alpha r_{j^*}^- + (1-\alpha) r_{j^*}^+,
\end{align}
\quad and
\begin{align}
\label{splitting_criterion_2}
    \#\{ i\in \cI:\bx_i\in L_l \} \geq \beta f(c), \quad l=1, 2.
\end{align}
\quad Here, $L_1 = \{\bx\in \R^d: \bx^{(j^*)} \leq z\}$,  $L_2 = \{\bx\in \R^d: \bx^{(j^*)} > z\}$, 
$\bar{Y}_1, \bar{Y}_2$ are average responses of $\S_\cJ$ within child nodes $L_1, L_2$.
\STATE \textbf{Otherwise:}
\STATE \quad Select the best split value $z^*$ and direction $j^*$ by optimizing the MSE:
\begin{align*}
    z^*, j^* = \arg\min_{z, j}\sum_{i \in \cJ} ( Y_i - \bar{Y}_1 )^2 \1 \{ \bx_i \in L_1 \} + \sum_{i \in \cJ} (Y_i - \bar{Y}_2 )^2 \1 \{ \bx_i \in L_2 \},
\end{align*}
\quad ensuring \eqref{splitting_criterion_1} and \eqref{splitting_criterion_2}. Here, $L_1 = \{\bx\in \R^d: \bx^{(j)} \leq z\}$,  $L_2 = \{\bx\in \R^d: \bx^{(j)} > z\}$, 
$\bar{Y}_1, \bar{Y}_2$ are average responses of $\S_\cJ$ within child nodes $L_1, L_2$.
\RETURN Split direction $j^*$, split value $z^*$, leaf nodes $L_1, L_2$.
\end{algorithmic}
\end{algorithm}

The setting of two disjoint datasets $\S_\cI$ and $\S_\cJ$ satisfies Definition \ref{def:honest}. Moreover, the probability $\kappa$ for randomly selecting a direction satisfies Definition \ref{def:random_split}, while constraint \eqref{splitting_criterion_1} satisfies Definition \ref{def:spatially_balanced} and constraint \eqref{splitting_criterion_2} satisfies Definition \ref{def:sample_balanced}.

Using this criterion, we can grow a Regular Tree from the given data. The tree growing procedure is presented in Algorithm \ref{algorithm:RegularTree}. 
This algorithm begins by dividing the dataset $\S$ into two disjoint datasets $\S_\cI$ and $\S_\cJ$, and recursively applies Algorithm \ref{algorithm:splitting_criterion} to split the space until each node with depth $c$ contains $\beta f(c-1)$ to $f(c)-1$ samples from $\S_\cI$.
It can be shown that the tree grown by Algorithm \ref{algorithm:RegularTree} satisfies all four definitions.

\begin{algorithm}[htbp]
\caption{Regular Tree}
\label{algorithm:RegularTree}
\renewcommand{\algorithmicrequire}{\textbf{Input:}}
\renewcommand{\algorithmicendfor}{\textbf{Until:}}
\renewcommand{\algorithmicfor}{\textbf{Repeat:}}
\begin{algorithmic}[1]
\REQUIRE Dataset $\S=(\bx_i, Y_i)_{i=1}^n$, spatially balanced parameter $\alpha$, random-split parameter $\kappa$, and sample balanced parameter $f(c), \beta$.

\STATE Divide $\S$ into $\S_\cI$ and $\S_\cJ$ with $|\cI| = \lfloor n / 2\rfloor$ and $|\cJ| = n - \lfloor n / 2\rfloor$. 

\FOR{each current node $L \subseteq [0, 1]^d$}
\STATE Use Algorithm\ref{algorithm:splitting_criterion} with $\S^L_\cI = \{(\bx_i, Y_i): i\in \cI, \bx_i \in L \}$, $\S^L_\cJ = \{(\bx_i, Y_i): i\in \cJ, \bx_i \in L \}$ to find the split direction $j^*$, split value $z^*$. Divide $L$ into $L_1, L_2$.
\ENDFOR \quad Each node with depth $c$ contains $\beta f(c-1)$ to $f(c)-1$ samples in $\S_\cI$.
\RETURN A Regular Tree.
\end{algorithmic}
\end{algorithm}

\section{Regular Tree Search}\label{sec:regular_tree_search}
The Regular Tree Search algorithm involves two stages. In Stage 1, design points are sampled uniformly from the entire feature space $\cX$, while in Stage 2, they are selected sequentially. The algorithm is outlined in Algorithm \ref{algorithm:Regular_Tree_Search}.

\begin{algorithm}[htbp]
\caption{Regular Tree Search}
\begin{algorithmic}[1]
\label{algorithm:Regular_Tree_Search}
\REQUIRE Total sample budget $N$, spatially balanced parameter $\alpha $, random-split parameter $\kappa$, sample balanced parameter $f(c), \beta$, and stage 1 sample budget $N_0<N$.
\STATE Uniformly sample $N_0$ points from the space $\cX$. Simulate at each point to get $\S$. Set $n = N_0$.
\STATE Divide $\S$ into $\S_\cI$ and $\S_\cJ$ with $|\cI| = \lfloor n / 2\rfloor$ and $|\cJ| = n - \lfloor n / 2 \rfloor$. Build up a Regular tree $T$ using Algorithm \ref{algorithm:RegularTree}. Let $\cL$ be the set of leaf nodes of $T$.
\WHILE{$n < N$}
\STATE Select a leaf node $\hat{L} \in \cL$ using a tree search strategy.
\STATE Uniformly samples a point $\bx_{n+1}$ from the space of leaf node $\hat{L}$, simulate at point $\bx_{n+1}$ to get $Y_{n+1}$. 
\STATE $\S_\cI = \S_\cI \cup (\bx_{n+1}, Y_{n+1})$.
\IF{$\# \{ i \in \cI, \bx_i \in \hat{L} \} = f(c)$, where $c$ is the depth of $\hat{L}$}
\STATE Uniformly sample $\max \{ f(c) - \# \{ j \in \cJ: \bx_j \in \hat{L} \}, 0\}$ points, simulate at each point once, set $n = n + \max \{ f(c) - \# \{ j \in \cJ: \bx_j \in \hat{L} \}, 0\}$, and add these points to $\S_\cJ$. 
\STATE Use Algorithm \ref{algorithm:splitting_criterion} with $\S^{\hat{L}}_\cI = \{(\bx_i, Y_i): i\in \cI, \bx_i \in \hat{L} \}$, $\S^{\hat{L}}_\cJ = \{(\bx_i, Y_i): i\in \cJ, \bx_i \in \hat{L} \}$ to find the split direction $j^*$ and split value $z^*$. Divide $\hat{L}$ into $L_1, L_2$.
\STATE $\cL = \cL \cup L_1 \cup L_2 \backslash \hat{L}$.
\ENDIF
\STATE $n = n + 1 $.
\ENDWHILE
\RETURN The leaf node $L^* \in \cL $ with the largest sample mean
\begin{align}
\label{equation:optimal_solution}
    L^* = \arg\max_{L\in \cL}\frac{\sum_{i\in \cI} \1{\{ \bx_i\in L \}} Y_i}{ \sum_{i\in \cI} \1{ \{ \bx_i\in L \} } },
\end{align}
and its corresponding optimal value estimate
\begin{align}
\label{equation:optimal_value}
    \bar{Y}_{L^*} = \frac{\sum_{i\in \cI} \1{\{ \bx_i\in L^* \}} Y_i}{ \sum_{i\in \cI} \1{ \{ \bx_i\in L^* \} } }.
\end{align}
\end{algorithmic}
\end{algorithm}

In Stage 1, we allocate an initial budget of $N_0$ function evaluations to uniformly sample design points from the entire feasible domain $\cX$. Each sampled point is simulated once, and the resulting dataset is split into two disjoint subsets, $\cS_\cI$ for estimating leaf values, and $\cS_\cJ$ for determining the split direction and value. A Regular Tree can be grown via Algorithm \ref{algorithm:RegularTree}.
Stage 1 serves as a warm-up phase, designed to provide an initial exploration of the space. This phase is analogous to the warm-up procedure in Bayesian optimization \cite[see, e.g.,][]{frazier2018tutorial}.
After this stage, we can construct an initial space partitioning of the feasible domain. 
Increasing the budget allocated to Stage 1 allows for a finer initial partitioning, providing a better global approximation. 
However, if the entire budget is exhausted during this phase, the algorithm degenerates into a regression algorithm, aiming to predict the function over the feasible domain rather than to identify optimal regions.

After Stage 1, the remaining sample budget is $N-N_0$. In each iteration, we select the next design point to be simulated using the tree search strategy. A tree search strategy is to select a leaf node $\hat{L} \in \cL$ while balancing exploration and exploitation. Once a leaf node $\hat{L}$ is selected, a new point is uniformly sampled from its region and evaluated, and the resulting observation is added to $\S_\cI$. 
If the number of $\cI$-samples in $\hat{L}$ reaches the threshold $f(c)$, where $c$ denotes the depth of the leaf $\hat{L}$, a refinement of the partition is triggered. We uniformly sample points in $\hat{L}$ and add to $\S_\cJ$ until the number of $\cJ$-samples in $\hat{L}$ reaches the threshold $f(c)$. Then we partition $\hat{L}$ into two new children, $L_1$ and $L_2$, using Algorithm \ref{algorithm:splitting_criterion}. The set of leaf nodes $\cL$ is then updated.

Once the budget is exhausted, the algorithm returns the leaf node $L^* \in \cL$ with the highest sample mean of responses in $\S_\cI$.

An effective tree search strategy for Stage 2 employs the UCT algorithm, which recursively selects child nodes from the root (i.e., the feasible domain $\cX$) until reaching a leaf node. UCT operates on a selection criterion analogous to the UCB principle in multi-armed bandit problems, which is a highly efficient approach widely adopted in tree search algorithms.
At a parent node $\hat{L}$ with children node $L_1$ and $L_2$, the UCB is computed for each child:
\begin{align*}
    \text{ucb}_l = \bar{Y}_{L_l} + C_p \sqrt{\frac{2 \log n_p}{n_l}}, \quad l=1, 2,
\end{align*}
where $\bar{Y}_{L_l}$ is the sample mean of responses in $L_l$ computed from $\S_\cI$, $n_l$ is the number of samples in $L_l$ within $\S_\cI$, and $n_p$ is the number of samples in the parent node $\hat{L}$ within $\S_\cI$. The constant $C_p$ is a tunable hyperparameter that controls the trade-off between exploration and exploitation. Setting $C_p=0$ results in a purely greedy search, which always selects the node with the largest sample mean. For the larger $C_p$, the algorithm reduces to a fully exploratory strategy that behaves like uniform sampling, eventually expanding the entire tree.
At each step, the child node with the larger UCB value is selected as the new $\hat{L}$, and the process continues until a leaf node is reached. 

\section{Convergence Analysis}\label{sec:convergence_analysis}
In this section, we establish the convergence properties of Algorithm \ref{algorithm:Regular_Tree_Search}. Our objective is to show that as the simulation budget $N$ increases, the algorithm identifies a region that concentrates around the true optimal point $\bx^*$, and the corresponding estimated optimal value converges to the global maximum $\mu(\bx^*)$. Specifically, we prove that:
\begin{align*}
    \dist(L_N^*, \bx^* ) \xrightarrow{P} 0,\\
    \bar{Y}_{L_N^*} \xrightarrow{P}  \mu(\bx^*),
\end{align*}
where $L_N^*$ denotes the zone (leaf node) selected by Algorithm \ref{algorithm:Regular_Tree_Search} with the largest sample mean when the sample budget is $N$, and $\dist(L_N^*, \bx^* ) = \inf_{\bx\in L_N^*} \| \bx - \bx^* \|$.

We begin with Lemma \ref{lemma:subgaussian}, which establishes that for any region $L$, the random variable $Y(L)$ is sub-Gaussian. This property plays a crucial role in our analysis, as it enables the use of concentration inequalities to bound the deviation between the sample mean $\bar{Y}_L$ and its expected value $\E[Y(L)]$.
\begin{lemma}
\label{lemma:subgaussian}
    Suppose that Assumptions \ref{assum:bound} and \ref{assum:local_increase} hold. For any hyperrectangle $L \subseteq \cX$, $Y(L)$ is sub-Gaussian.
\end{lemma}

Theorem \ref{theorem:sol_gap_stochastic_subgaussian} states that the selected zone $L_N^*$ asymptotically shrinks toward the true optimal point under some conditions.
\begin{theorem}
\label{theorem:sol_gap_stochastic_subgaussian}
     Suppose that Assumptions \ref{assum:bound}, \ref{assum:local_increase}, and \ref{assum:local_decay} hold. Consider the tree search strategy that satisfies the following conditions:
    \begin{itemize}
        \item The depth of each leaf is bounded between $h^-(N)$ and $h^+(N)$.
        \item $\lim_{N \to \infty} h^-(N) = \infty$.
        \item $\lim_{N \to \infty}\frac{h^+(N)}{f(h^-(N))} = 0$.
    \end{itemize}
    Then, $\dist(L_N^*, \bx^*)$ converges in probability to 0.
\end{theorem}

Theorem \ref{theorem:sol_gap_stochastic_subgaussian} contains three conditions:
\begin{enumerate}
    \item During the algorithm, the depth of any leaf node lies between $h^-(N)$ and $h^+(N)$.
    \item The tree search strategy is designed so that $h^-(N)\to \infty$ as the total budget $N \to \infty$, ensuring that the partition is refined ever deeper over time. 
    \item The ratio $\frac{h^+(N)}{f(h^-(N))} \to 0$ as $N \to \infty$ imposes a balance between the maximal tree depth $h^+(N)$ and the number of samples required for splitting. Intuitively, this condition prevents overfitting due to overly deep trees without sufficient data, thereby maintaining a balance between partition refinement and statistical reliability.
\end{enumerate}

These three conditions imply that the sample balanced function $f(c)$, which determines the number of samples required to split a node at depth $c$, diverges as $c \to \infty$. This ensures deeper leaf nodes receive enough observations to yield statistically reliable estimates.

The exploration-exploitation trade-off in the UCT algorithm ensures that each node is selected infinitely many times, implying that $\lim_{N \to \infty} h^-(N) = \infty$. In practice, if $h^-(N)$ and $h^+(N)$ are of the same order, then setting $f(c)$ to be the order of $c \log c$ satisfies the conditions of Theorem \ref{theorem:sol_gap_stochastic_subgaussian}. However, designing a tree search strategy that explicitly guarantees the desired behaviors of $h^-(N)$ and $h^+(N)$ remains an open problem and a promising direction for future research.

Building on Theorem \ref{theorem:sol_gap_stochastic_subgaussian}, Theorem \ref{theorem:optimal_gap_stochastic} establishes that the sample mean of the selected region, $\bar{Y}_{L_N^*}$, converges in probability to the true optimal value $\mu(\bx^*)$.
\begin{theorem}
\label{theorem:optimal_gap_stochastic}
    Under the conditions of Theorem \ref{theorem:sol_gap_stochastic_subgaussian}, $\bar{Y}_{L_N^*} - \mu(\bx^*)$ converges in probability to $0$.
\end{theorem}

\section{Numerical Experiments}\label{sec:numerical_experiments}
In this section, we conduct numerical experiments to evaluate the performance of the Regular Tree Search algorithm with UCT.
We compare our algorithm with three global convergence algorithms: the ASR algorithm from \cite{andradottir2010adaptive}, and the IHR-SO and AP-SO algorithms from \cite{kiatsupaibulSingleObservationAdaptive2018}.

Consider the Rastrigin function:
\begin{align*}
    \mu(\bx) = 10d + \sum_{i=1}^d [x_i^2 - 10 \cos (2\pi x_i)],
\end{align*}
where $\bx = (x_1, \dots, x_d)$.
The feasible domain is $\cX = [-5, 5]^d$. The objective is to find the global minimum at $\bx^*=(0, \dots, 0)$ with $\mu(\bx^*)=0$. 
This function has multiple local minima. Let $\epsilon(\bx) \sim N(0, 1)$. For illustration purposes, we consider a problem formulation that differs from our algorithm design by converting the original maximization problem into a minimization problem. This change is equivalent to solving the problem of maximizing the negative of the Rastrigin function.


The function is tested for different dimensions ($d=2, 5, 10$) to evaluate the algorithm's performance. 
The total sample budget $N$ is set to $500d$. 
For Regular Tree Search, the first stage sample budget is set to $N_0=150d$. Smaller value of parameters $\alpha$, $\kappa$ and $\beta$, combined with larger $C_p$ promote exploration. We adopt the default settings: the spatially balanced parameter to $\alpha=0.1$, the random-split parameter to $\kappa=0.1$, the UCT parameter to $C_p=2$, splitting parameter to $\beta=1/3$. 
Theorem \ref{theorem:optimal_gap_stochastic} shows that one may choose $f(c)$ on the order of $ c\log c$, although theoretical guarantees for the order of $h^-(N)$ and $h^+(N)$ remain incomplete. To balance finite-sample performance against the risk of premature overpartitioning in early iterations, we define the default sample balanced parameter $f(c) = \max(c\log c, 15)$.
The ASR, IHR-SO, and AP-SO algorithms use the same parameter settings as those in \cite{andradottir2010adaptive} and \cite{kiatsupaibulSingleObservationAdaptive2018}. 

Each experiment is repeated 100 times. For each algorithm, we record the \textit{true function value of the optimal solution estimates} and the \textit{optimal value estimates}. 
Since Regular Tree Search ultimately selects a region, as shown in \eqref{equation:optimal_solution}, we use the midpoint of the selected region as the optimal solution estimate to facilitate comparison. 
The corresponding estimated optimal value is computed as in \eqref{equation:optimal_value}.
The true function value at the estimated optimal solution and the estimated optimal value for ASR, IHR-SO, and AP-SO algorithms are computed in the same manner as described in \cite{kiatsupaibulSingleObservationAdaptive2018}.
We report summary statistics over 100 replications, including the mean, Root MSE (RMSE), the best value, 25th, 50th, and 75th percentiles, as well as the worst value.
Table \ref{tab:comparison_true} presents the performance statistics of the true function value of the optimal solution estimates of four algorithms across different dimensions $d$, while Table \ref{tab:comparison_estimate} presents the corresponding performance statistics of the optimal value estimates.

\begin{table}[t]
    \centering
    \caption{Statistics of the true function value of the optimal solution estimate of four algorithms.}
    \label{tab:comparison_true}
    \begin{tabular}{c c c c c c c c c}
        \toprule
        & & & & & \multicolumn{3}{c}{Quantile of value} &  \\
        \cmidrule(lr){6-8}
        $d$ & Algorithm & Mean & RMSE & Best & 25\% & 50\% & 75\% & Worst \\
        \midrule
        \multirow{4}{*}{2} & Regular Tree Search & \multicolumn{1}{c}{3.20} & \multicolumn{1}{c}{4.17} & \multicolumn{1}{c}{0.03} & \multicolumn{1}{c}{1.28} & \multicolumn{1}{c}{2.11} & \multicolumn{1}{c}{4.99} & \multicolumn{1}{c}{10.23} \\
         & ASR & \multicolumn{1}{c}{5.23} & \multicolumn{1}{c}{6.16} & \multicolumn{1}{c}{0.31} & \multicolumn{1}{c}{3.02} & \multicolumn{1}{c}{4.78} & \multicolumn{1}{c}{6.73} & \multicolumn{1}{c}{20.50} \\
         & IHR-SO & \multicolumn{1}{c}{36.54} & \multicolumn{1}{c}{36.74} & \multicolumn{1}{c}{22.61} & \multicolumn{1}{c}{34.24} & \multicolumn{1}{c}{37.18} & \multicolumn{1}{c}{39.69} & \multicolumn{1}{c}{43.80} \\
         & AP-SO & \multicolumn{1}{c}{35.77} & \multicolumn{1}{c}{36.16} & \multicolumn{1}{c}{21.06} & \multicolumn{1}{c}{33.52} & \multicolumn{1}{c}{37.06} & \multicolumn{1}{c}{39.77} & \multicolumn{1}{c}{43.83} \\
         \midrule
        \multirow{4}{*}{5} & Regular Tree Search & \multicolumn{1}{c}{9.16} & \multicolumn{1}{c}{10.08} & \multicolumn{1}{c}{1.82} & \multicolumn{1}{c}{6.10} & \multicolumn{1}{c}{8.77} & \multicolumn{1}{c}{12.22} & \multicolumn{1}{c}{20.15} \\
         & ASR & \multicolumn{1}{c}{26.63} & \multicolumn{1}{c}{28.16} & \multicolumn{1}{c}{9.80} & \multicolumn{1}{c}{20.37} & \multicolumn{1}{c}{25.42} & \multicolumn{1}{c}{31.36} & \multicolumn{1}{c}{52.27} \\
         & IHR-SO & \multicolumn{1}{c}{45.79} & \multicolumn{1}{c}{47.29} & \multicolumn{1}{c}{17.41} & \multicolumn{1}{c}{39.17} & \multicolumn{1}{c}{43.91} & \multicolumn{1}{c}{51.62} & \multicolumn{1}{c}{82.13} \\
         & AP-SO & \multicolumn{1}{c}{59.87} & \multicolumn{1}{c}{60.92} & \multicolumn{1}{c}{19.98} & \multicolumn{1}{c}{52.52} & \multicolumn{1}{c}{60.30} & \multicolumn{1}{c}{67.56} & \multicolumn{1}{c}{87.37} \\
         \midrule
        \multirow{4}{*}{10} & Regular Tree Search & \multicolumn{1}{c}{27.86} & \multicolumn{1}{c}{30.08} & \multicolumn{1}{c}{9.58} & \multicolumn{1}{c}{19.02} & \multicolumn{1}{c}{27.21} & \multicolumn{1}{c}{35.98} & \multicolumn{1}{c}{65.51} \\
         & ASR & \multicolumn{1}{c}{76.68} & \multicolumn{1}{c}{79.04} & \multicolumn{1}{c}{39.78} & \multicolumn{1}{c}{62.34} & \multicolumn{1}{c}{73.60} & \multicolumn{1}{c}{88.46} & \multicolumn{1}{c}{134.56} \\
         & IHR-SO & \multicolumn{1}{c}{92.28} & \multicolumn{1}{c}{93.72} & \multicolumn{1}{c}{43.84} & \multicolumn{1}{c}{81.50} & \multicolumn{1}{c}{91.57} & \multicolumn{1}{c}{103.33} & \multicolumn{1}{c}{138.44} \\
         & AP-SO & \multicolumn{1}{c}{83.55} & \multicolumn{1}{c}{86.23} & \multicolumn{1}{c}{39.01} & \multicolumn{1}{c}{68.23} & \multicolumn{1}{c}{82.47} & \multicolumn{1}{c}{96.12} & \multicolumn{1}{c}{137.12} \\
        \bottomrule
    \end{tabular}
\end{table}

\begin{table}[t]
    \centering
    \caption{Statistics of optimal value estimates of four algorithms.}
    \label{tab:comparison_estimate}
    \begin{tabular}{c c c c c c c c c}
        \toprule
        & & & & & \multicolumn{3}{c}{Quantile of value} &  \\
        \cmidrule(lr){6-8}
        $d$ & Algorithm & Mean & RMSE & Best & 25\% & 50\% & 75\% & Worst \\
        \midrule
        \multirow{4}{*}{2} & Regular Tree Search & \multicolumn{1}{c}{2.72} & \multicolumn{1}{c}{3.79} & \multicolumn{1}{c}{0.02} & \multicolumn{1}{c}{0.68} & \multicolumn{1}{c}{1.69} & \multicolumn{1}{c}{4.29} & \multicolumn{1}{c}{9.78} \\
         & ASR & \multicolumn{1}{c}{5.23} & \multicolumn{1}{c}{6.17} & \multicolumn{1}{c}{0.34} & \multicolumn{1}{c}{4.01} & \multicolumn{1}{c}{4.85} & \multicolumn{1}{c}{6.74} & \multicolumn{1}{c}{20.64} \\
         & IHR-SO & \multicolumn{1}{c}{19.70} & \multicolumn{1}{c}{19.71} & \multicolumn{1}{c}{17.89} & \multicolumn{1}{c}{19.16} & \multicolumn{1}{c}{19.76} & \multicolumn{1}{c}{20.25} & \multicolumn{1}{c}{21.56} \\
         & AP-SO & \multicolumn{1}{c}{21.43} & \multicolumn{1}{c}{21.44} & \multicolumn{1}{c}{18.00} & \multicolumn{1}{c}{20.99} & \multicolumn{1}{c}{21.39} & \multicolumn{1}{c}{21.98} & \multicolumn{1}{c}{23.89} \\
         \midrule
        \multirow{4}{*}{5} & Regular Tree Search & \multicolumn{1}{c}{8.97} & \multicolumn{1}{c}{9.91} & \multicolumn{1}{c}{1.29} & \multicolumn{1}{c}{6.07} & \multicolumn{1}{c}{8.55} & \multicolumn{1}{c}{11.96} & \multicolumn{1}{c}{20.22} \\
         & ASR & \multicolumn{1}{c}{26.60} & \multicolumn{1}{c}{28.13} & \multicolumn{1}{c}{9.82} & \multicolumn{1}{c}{20.37} & \multicolumn{1}{c}{25.40} & \multicolumn{1}{c}{31.30} & \multicolumn{1}{c}{52.19} \\
         & IHR-SO & \multicolumn{1}{c}{41.51} & \multicolumn{1}{c}{41.55} & \multicolumn{1}{c}{37.66} & \multicolumn{1}{c}{40.26} & \multicolumn{1}{c}{41.51} & \multicolumn{1}{c}{42.79} & \multicolumn{1}{c}{46.40} \\
         & AP-SO & \multicolumn{1}{c}{48.88} & \multicolumn{1}{c}{48.92} & \multicolumn{1}{c}{45.17} & \multicolumn{1}{c}{47.09} & \multicolumn{1}{c}{48.94} & \multicolumn{1}{c}{50.38} & \multicolumn{1}{c}{55.39} \\
         \midrule
        \multirow{4}{*}{10} & Regular Tree Search & \multicolumn{1}{c}{31.64} & \multicolumn{1}{c}{33.12} & \multicolumn{1}{c}{12.87} & \multicolumn{1}{c}{23.77} & \multicolumn{1}{c}{31.14} & \multicolumn{1}{c}{38.89} & \multicolumn{1}{c}{59.61} \\
         & ASR & \multicolumn{1}{c}{76.62} & \multicolumn{1}{c}{78.98} & \multicolumn{1}{c}{39.57} & \multicolumn{1}{c}{62.37} & \multicolumn{1}{c}{73.72} & \multicolumn{1}{c}{88.32} & \multicolumn{1}{c}{134.47} \\
         & IHR-SO & \multicolumn{1}{c}{84.83} & \multicolumn{1}{c}{85.02} & \multicolumn{1}{c}{73.59} & \multicolumn{1}{c}{80.94} & \multicolumn{1}{c}{84.19} & \multicolumn{1}{c}{88.57} & \multicolumn{1}{c}{101.80} \\
         & AP-SO & \multicolumn{1}{c}{93.95} & \multicolumn{1}{c}{94.18} & \multicolumn{1}{c}{80.46} & \multicolumn{1}{c}{88.89} & \multicolumn{1}{c}{94.15} & \multicolumn{1}{c}{97.91} & \multicolumn{1}{c}{106.99} \\
        \bottomrule
    \end{tabular}
\end{table}

For the $d=2$ case, Regular Tree Search outperforms the competing methods. The algorithm achieves a mean true function value of 3.20 with an RMSE of 4.17, and its best-case performance is exceptionally low at 0.03. The quantile values (1.28 at the 25th percentile, 2.11 at the 50th, and 4.99 at the 75th) exhibit a narrow spread, indicating that Regular Tree Search not only delivers high accuracy but also maintains robust and consistent performance. In comparison, ASR records a higher mean of 5.23 and RMSE of 6.16, while both IHR-SO and AP-SO yield substantially larger mean values (36.54 and 35.77, respectively).
A similar trend is observed for the optimal value estimates. Here, Regular Tree Search attains a mean of 2.72 and an RMSE of 3.79 with a best-case value of 0.02, reinforcing its superior accuracy and consistency in low-dimensional settings.

When the dimensionality increases to $d=5$ and $d=10$, Regular Tree Search remains competitive despite the challenges of higher dimensions. For $d=5$, it has a true function value mean of 9.16 and an RMSE of 10.08, with quantiles expanding from 6.10 at the 25th to 12.22 at the 75th. At $d=10$, it maintains a mean true function value of 27.86 and an RMSE of 30.08, with quantile values ranging from 19.02 (25th) to 35.98 (75th). In both cases, the competing algorithms experience sharper performance deterioration (with ASR, IHR-SO, and AP-SO recording mean values of 26.63, 45.79, 59.87 at $d=5$ and 76.68, 92.28, 83.55 at $d=10$, respectively), while Regular Tree Search continues to deliver lower error levels and more robust performance.
Similar trends are confirmed by the optimal value estimates, where Regular Tree Search maintains a lower mean and error level relative to its counterparts.

When comparing the true function value of the optimal solution estimate with the optimal value estimates, a clear pattern emerges: both Regular Tree Search and SO algorithms exhibit noticeable gaps between their estimated and actual values, whereas ASR shows virtually no discrepancy. This difference arises because ASR performs multiple independent simulations at each design point and directly computes its sample mean, while Regular Tree Search, IHR-SO, and AP-SO each run only one simulation per point and then use different strategies to estimate the function value.

In particular, for $d = 2$ and $d = 5$, IHR-SO and AP-SO tend to underestimate the true function value. This is likely because these algorithms estimate values by averaging over a neighborhood, which can weaken the impact of the true optimum, especially if the function has sharp peaks. If the initial neighborhood size is too large, the estimate may become inaccurate.

Regular Tree Search enforces an honest sampling rule, ensuring that within‐leaf sample means are unbiased estimators of the region's average. However, this does not guarantee an unbiased estimate at the center of the region. Consequently, in the $d=2$ scenario, Regular Tree Search's estimated optimal value falls below the true optimal value for two reasons: first, we report the function value at the region's center rather than the region's actual mean; second, because Regular Tree Search always selects the leaf with the lowest estimated mean, this selection step introduces a negative selection bias in the reported optimal value estimate.
As dimensionality increases to $d=5$ and $d=10$, the relationship between estimated and true values becomes less predictable. Sometimes, the estimates exceed the true values; sometimes, they fall short. This increased variability can be attributed to the coarser partitioning required in higher dimensions, which results in larger regions. When a region is large, the discrepancy between its center and its true mean is also large, leading to estimation errors.

\section{Conclusion}\label{sec:conclusion}
This paper introduces Regular Tree Search, a class of random search algorithms for non-convex simulation optimization that integrates adaptive sampling with adaptive space partitioning. The algorithm initiates with a regular tree structure constructed from uniformly sampled points and iteratively selects leaf nodes by dynamically balancing exploration and exploitation through a tree search strategy. A new sample is drawn uniformly within the selected leaf, and when the number of samples in a node exceeds a threshold determined by its depth, an adaptive splitting criterion triggers further partitioning, progressively focusing the search on promising regions. The algorithm ultimately selects the best empirical solution from the leaves when the sampling budget is exhausted.
To guide the search process, the framework incorporates the UCT algorithm. From a theoretical view, we establish global convergence under sub-Gaussian noise assumptions, relying on conditions concerning the optimality gap rather than global continuity. Numerical experiments validate the effectiveness of the algorithm, demonstrating its ability to accurately locate optimal solutions and estimate their values in complex, noisy optimization settings.

\bibliographystyle{plainnat}
\bibliography{mybibfile}

\newpage

\appendix
\section{Proof of Theorems and Lemmas}

\subsection{Proof of Lemma \ref{lemma:diam_depth}}

We first introduce Lemma \ref{lemma:cor_depth} and then prove Lemma \ref{lemma:diam_depth}.
\begin{lemma}
\label{lemma:cor_depth}
    Let $T$ be a random-split, $\alpha$-spatially balanced tree, and let $L$ be a leaf of $T$ with depth $c$. Suppose that leaf $L = \bigotimes_{j=1}^d \left[ r_j^-, r_j^+ \right]$, where $0 \leq r_j^- < r_j^+ \leq 1$. For each coordinate $j$,  define the diameter of $L$ along the $j$-th coordinate as $\diam_j(L) = r_j^+ - r_j^-$. Then, for any $\lambda \in (0, 1)$:
    \begin{align*}
        \P \left[\diam_j(L) \geq (1-\alpha)^{(1-\lambda)\frac{\kappa}{d}c}\right] \leq \exp \left\{-\frac{\lambda^2}{2} \frac{\kappa}{d} c \right\}.
    \end{align*}
\end{lemma}
\begin{proof}
    Each time a split is performed on a node, the depth of its child node increases by $1$. Therefore, the depth of a leaf $L$ corresponds to the total number of splits that have occurred along the path to $L$. Therefore, depth $c$ is the number of splits of the leaf $L$. Let $c_j$ be the number of splits at the leaf $L$ along the $j$-th coordinate. Since the tree is constructed as a random-split tree, where each split selects a coordinate bounded below by $\frac{\kappa}{d}$, and $\kappa$ represents the parameter for the random-split, the random variable $c_j$ is stochastically bounded below by a Binomial distribution:
    \begin{align*}
        c_j \stackrel{d}{\geq} \text{Binom} \left( c ; \frac{\kappa}{d}\right),
    \end{align*}
    where $\stackrel{d}{\geq}$ denotes that random variable $c_j$ is stochastically greater than or equal to the binomial random variable. Specifically, this means:
    \begin{align}
    \label{ineq:sto_greater}
        \P( c_j \geq t ) \geq \P \left( \text{Binom} \left( c ; \frac{\kappa}{d}\right) \geq t \right), \quad \forall t \in \R.
    \end{align}

    Next, we apply Chernoff’s inequality to bound the tail probability of $c_j$. Chernoff’s inequality for a binomial random variable $X \sim \text{Binom}(n, p)$  states that for any $\lambda \in (0, 1)$,
    \begin{align*}
        \P \left[ X \leq (1 - \lambda) \E[X] \right] \leq \exp \left\{ - \frac{\lambda^2}{2} \E[X] \right\}.
    \end{align*}
    
    In our case, since $c_j$ is stochastically greater than or equal to $\text{Binom}(c, \frac{\kappa}{d})$, we can apply Chernoff’s inequality to the binomial random variable. From \eqref{ineq:sto_greater}, for any $\lambda > 0$,
    \begin{align}
    \label{ineq:lemma1_chernoff}
         \P \left[ c_j \leq (1 - \lambda) \frac{\kappa}{d} c \right] \leq \P \left[ \text{Binom} \left( c ; \frac{\kappa}{d}\right) \leq (1 - \lambda) \frac{\kappa}{d} c \right] \leq \exp \left\{ - \frac{\lambda^2}{2} \frac{\kappa}{d} c \right\}.
    \end{align}

    From the definition of $\alpha$-spatially balanced, each split along the $j$-th coordinate reduces the diameter of $L$ in that coordinate by at least a factor of $ 1-\alpha$. Hence, the diameter of $L$ along the $j$-th coordinate satisfies:
    \begin{align*}
        \diam_j(L) \leq \left( 1 - \alpha \right)^{ c_j }.
    \end{align*}
    
    Combining the result with the above Chernoff bound \eqref{ineq:lemma1_chernoff}, we obtain:
    \begin{align*}
        \P\left[\diam_j(L) \geq (1-\alpha)^{(1-\lambda)\frac{\kappa}{d}c}\right] &\leq \P\left[c_j \leq (1 - \lambda) \frac{\kappa}{d} c \right] \leq \exp\left\{-\frac{\lambda^2}{2} \frac{\kappa}{d} c \right\}.
    \end{align*}
\end{proof}

\textbf{Proof of Lemma \ref{lemma:diam_depth}.}
\begin{proof}
    To extend bound in Lemma \ref{lemma:cor_depth} to the 2-norm of $\diam(L)$, note that:
    \begin{align*}
        \diam (L) = \sqrt{\sum_{j=1}^d \diam_j^2(L)}.
    \end{align*}

    Using a union bound over all $d$ coordinates, we have:
    \begin{align*}
        \P \left[\diam (L) \geq \sqrt{d} (1-\alpha)^{(1-\lambda) \frac{\kappa}{d} c}\right] &\leq \P \left[\bigcup_{j=1}^d \left\{ \diam_j (L) \geq (1-\alpha)^{(1-\lambda) \frac{\kappa}{d} c} \right\}\right] \\
        &\leq \sum_{j=1}^d \P \left[\diam_j (L) \geq (1-\alpha)^{(1-\lambda) \frac{\kappa}{d} c}\right] \\
        &\leq d \exp \left\{-\frac{\lambda^2}{2} \frac{\kappa}{d} c \right\}.
    \end{align*}
\end{proof}

\subsection{Proof of Lemma \ref{lemma:subgaussian}}
\begin{proof}
    From Assumption \ref{assum:local_increase},
    \begin{align*}
        \mu(\bx^*) - C \cdot \| \bx - \bx^* \| \leq \mu(\bx) \leq \mu(\bx^*).
    \end{align*}

    Define $C_\mu = \max \{ |\mu(\bx^*) - C|, |\mu(x^*)| \}$, then it follows that $|\mu(\bx)| \leq C_\mu$ for all $\bx \in L$. 
    To prove $Y(L)$ is sub-Gaussian, we analyze its moment-generating function.
    For any $\lambda \in \mathbb{R}$:
    \begin{align*}
        \mathbb{E}\left[ e^{\lambda (Y(L) - \E[Y(L)]) } \right] &= \mathbb{E}_{\bX \sim \text{U}(L)}\left[ \mathbb{E}_{Y|\bX}\left[ e^{\lambda Y} \mid \bX \right] \right] \cdot e^{-\lambda \E[Y(L)]} \\
        &= \mathbb{E}_{\bX \sim \text{U}(L)}\left[ e^{\lambda \mu(\bX)} \cdot \mathbb{E}_{Y|\bX}\left[ e^{\lambda \epsilon(\bX)} \mid \bX \right] \right] \cdot e^{-\lambda \E[Y(L)]}.
    \end{align*}
    
    By the sub-Gaussianity of $\epsilon(\bx)$:
    \begin{align}
        \mathbb{E}_{Y|\bX}\left[ e^{\lambda \epsilon(\bX)} \mid \bX \right] \leq e^{\lambda^2 C_{\text{sg}}^2 / 2}.
    \end{align}
    
    Thus:
    \begin{align*}
        \E\left[ e^{\lambda (Y(L) - \E[Y(L)])} \right] &\leq \E_{\bX \sim \text{U}(L)}\left[ e^{\lambda \mu(\bX)} \right] \cdot e^{-\lambda \E[Y(L)]} \cdot e^{\lambda^2 C_{\text{sg}}^2 / 2} \\
        &= \E_{\bX \sim \text{U}(L)}\left[ e^{\lambda \mu(\bX)} \right] \cdot e^{-\lambda \E[\mu(\bX)]} \cdot e^{\lambda^2 C_{\text{sg}}^2 / 2}\\
        &= \E_{\bX \sim \text{U}(L)}\left[ e^{\lambda (\mu(\bX)-\E[\mu(\bX)])} \right] \cdot e^{\lambda^2 C_{\text{sg}}^2 / 2}.
    \end{align*}
    
    Since $|\mu(\bx)| \leq C_\mu$ for all $\bx \in L$, the random variable $\mu(\bX)$ is bounded. By the Hoeffding lemma, any bounded random variable \( Z \in [a, b] \) is sub-Gaussian with variance proxy $ \frac{(b - a)^2}{4} $. Here, $\mu(\bX) \in [-C_\mu, C_\mu]$, so it is sub-Gaussian with variance proxy $ C_\mu^2 $:
    \begin{align*}
        \mathbb{E}\left[ e^{\lambda (\mu(\bX) - \E[\mu(\bX)])} \right] \leq e^{\lambda^2 C_\mu^2 / 2}.
    \end{align*}
    
    Then,
    \begin{align*}
        \E\left[ e^{\lambda (Y(L) - \E[Y(L)])} \right] \leq e^{\lambda^2 (C_\mu^2 + C_{\text{sg}}^2) / 2}.
    \end{align*}
\end{proof}

\subsection{Proof of Theorem \ref{theorem:sol_gap_stochastic_subgaussian}}
\begin{proof}
Consider any $\epsilon > 0$. We want to show that $\P(\dist(L_N^*, \bx^*)>\epsilon)$ converge to 0.
    By Assumption \ref{assum:local_decay}, there exists $\delta(\epsilon) > 0$ such that for any $\bx$ satisfying $\| \bx - \bx^*\| > \epsilon$, we have:
    \begin{align}
    \label{ineq:epsilon_delta}
        \mu(\bx) < \mu(\bx^*) - \delta(\epsilon).
    \end{align}

    Let $\cL_N$ denote the set of leaf nodes when the total sample size is $N$, and let $L_N(\bx^*)$ denote the leaf containing $\bx^*$ when the total sample size is $N$.

    Let us analyze the probability of the selected zone being far from the optimal point:
    \begin{align}
        \P(\dist (L_N^*, \bx^*) > \epsilon) &= \P(\text{the distance of zone with maximum sample average and $\bx^*$ is larger than $\epsilon$})\nonumber\\
        &= \P(\exists L \in \cL_N: \{\dist(L, \bx^*) > \epsilon \} \cap \{ \bar{Y}_L > \bar{Y}_{L_N(\bx^*)} \}). \label{ineq:prob_exist_L}
    \end{align}
    
    For any $\eta > 0$, we can decompose \eqref{ineq:prob_exist_L}:
    \begin{align}
        &\P(\dist(L_N^*, \bx^*) > \epsilon) \nonumber\\
        &= \P(\exists L \in \cL_N: \{\dist(L, \bx^*) > \epsilon \} \cap \{ \bar{Y}_L > \bar{Y}_{L_N(\bx^*)}\} | \bar{Y}_{L_N(\bx^*)} < \mu(\bx^*) - \eta) \P(\bar{Y}_{L_N(\bx^*)} < \mu(\bx^*) - \eta) \nonumber\\
        &\quad + \P(\exists L \in \cL_N: \{\dist(L, \bx^*) > \epsilon\} \cap \{\bar{Y}_L > \bar{Y}_{L_N(\bx^*)}\} | \bar{Y}_{L_N(\bx^*)} \geq \mu(\bx^*) - \eta) \P(\bar{Y}_{L_N(\bx^*)} \geq \mu(\bx^*) - \eta) \nonumber\\
        &\leq \P(\bar{Y}_{L_N(\bx^*)} < \mu(\bx^*) - \eta) + \P(\exists L \in \cL_N: \{\dist(L, \bx^*) > \epsilon \} \cap \{ \bar{Y}_L \geq \mu(\bx^*) - \eta\})\nonumber\\
        &= \P( \mu(\bx^*) - \bar{Y}_{L_N(\bx^*)} > \eta) + \P(\cup_{L \in \cL_N} \{ \{ \dist(L, \bx^*) > \epsilon \} \cap \{ \bar{Y}_L \geq \mu(\bx^*) - \eta \} \})\nonumber\\
        &\leq \P \left(\mu(\bx^*) - \E \left[Y( L_N (\bx^*) ) |L_N(\bx^*)\right] > \eta/2 \right) + \P(| \E[Y( L_N (\bx^*) ) |L_N(\bx^*) ] - \bar{Y}_{L_N(\bx^*)} | > \eta/2  )  \nonumber\\
        &\quad + \P(\cup_{L \in \cL_N} \{ \{ \dist(L, \bx^*) > \epsilon \} \cap \{ \bar{Y}_L \geq \mu(\bx^*) - \eta \} \}), \label{ineq:theorem_1}
    \end{align}

    \textbf{The first term.}
    
    Note that the randomness in the probability is reflected in $L_N(\bx^*)$. From the law of total expectation, we have:
    \begin{align*}
        \E\left[Y (L) | L \right]  = \E\left[ \E\left[Y | \bX \sim \text{U}(L) \right]\right] = \E_{\bX \sim \text{U}(L)} \left[ \mu(\bX) \right].
    \end{align*}
    By Assumption \ref{assum:local_increase} and the fact that $\bx^* \in L_N(\bx^*)$, we have:
    \begin{align}
        \mu(\bx^*) - \E[Y( L_N (\bx^*) ) |L_N(\bx^*) ] &= \E [\mu(\bx^*) - \mu(\bX) | \bX \sim \text{U} (L_N(\bx^*))] \nonumber\\
        &\leq \E[C \cdot \dist(\bx^*, \bX) | \bX \sim \text{U} (L_N(\bx^*))] \nonumber \\
        &\leq C \cdot \diam(L_N(\bx^*)).    \label{ineq:local_smoothness}
    \end{align}

    Let $\tau_N = \sqrt{d} (1-\alpha)^{(1-\lambda) \frac{\kappa}{d} h^-(N)}$ denote the diameter threshold. We decompose the target probability into two cases:
    \begin{align}
        &\P \left(| \mu(\bx^*) - \E \left[ Y( L_N (\bx^*) ) | L_N(\bx^*)\right]  | > \eta/2 \right) \nonumber \\
        &\leq \P \left( \left(\mu(\bx^*) - \E[Y( L_N (\bx^*) ) |L_N(\bx^*)] \right) \cdot \1 \left\{ \diam(L_N(\bx^*))>\tau_N \right\} > \eta/4 \right) \nonumber\\
        &\quad + \P \left( \left(\mu(\bx^*) - \E[Y( L_N (\bx^*) ) |L_N(\bx^*)] \right) \cdot \1 \left\{ \diam(L_N(\bx^*)) \leq \tau_N \right\} > \eta/4 \right).\nonumber
    \end{align}

    For the first case:
    \begin{align}
        &\P \left( \left(\mu(\bx^*) - \E[Y( L_N (\bx^*) ) |L_N(\bx^*)] \right) \cdot \1 \{ \diam(L_N(\bx^*))> \tau_N \} > \eta/4 \right)\nonumber \\
        &= \P \left( \mu(\bx^*) - \E[Y( L_N (\bx^*) ) |L_N(\bx^*)] > \eta/4, \diam(L_N(\bx^*))> \tau_N \right)\nonumber \\
        &\leq \P \left( \diam(L_N(\bx^*))>\tau_N \right)\nonumber \\
        &= \P \left( \diam(L_N(\bx^*))>\sqrt{d} (1-\alpha)^{(1-\lambda) \frac{\kappa}{d}  h^-(N)} \right)\nonumber \\
        &\leq d \exp \left\{ -\frac{\lambda^2}{2}\frac{\kappa}{d} h^-(N) \right\}, \label{ineq:deterministic_first_first_case}
    \end{align}
    where the last inequality follows from Lemma \ref{lemma:diam_depth} and the minimum leaf depth $h^-(N)$.

    For the second case:
    \begin{align}
        &\P \left((\mu(\bx^*) - \E[Y( L_N (\bx^*) ) |L_N(\bx^*)]) \cdot \1 \left\{ \diam(L_N(\bx^*)) \leq \tau_N \right\} > \eta/4 \right) \nonumber\\
        &= \P \left( \mu(\bx^*) - \E[Y( L_N (\bx^*) ) |L_N(\bx^*)]>\eta/4, \diam(L_N(\bx^*)) \leq \tau_N  \right) \nonumber \\
        &\leq \P \left( C \cdot \diam(L_N(\bx^*)) > \eta/4, \diam(L_N(\bx^*)) \leq \tau_N \right), \nonumber
    \end{align}
    where the inequality follows from \eqref{ineq:local_smoothness}.
    For $N$ large enough, $\tau_N < \eta/4C$, making this probability term vanish.

    Combining both cases, for $N$ large enough,
    \begin{align}
        \text{the first term} \leq d \exp \left\{ -\frac{\lambda^2}{2}\frac{\kappa}{d} h^-(N) \right\}. \label{ineq:first_term_combin}
    \end{align}
    
    RHS of \eqref{ineq:first_term_combin} converges to 0 as $N\to \infty$ because of $\lim_{N \to \infty} h^-(N) = \infty$. Then, the first term of \eqref{ineq:theorem_1} converges to $0$.

    \textbf{The second term.}

    Note that the randomness in the probability is reflected in $L_N(\bx^*)$ and $\bar{Y}$.
    From the Lemma \ref{lemma:subgaussian}, we know that for any $L$, $Y(L)$ is subgaussian with variance proxy $C_\mu^2 + C_{\text{sg}}^2$. 
    From the Hoeffding inequality and $L_N(\bx^*)$ contains at least $\beta f(h^-(N)-1)$ samples, we have:
    \begin{align}
    \label{ineq:subgaussian_term2}
        &\P\left( \left| \E[Y( L_N (\bx^*) ) | L_N(\bx^*) ] - \bar{Y}_{L_N(\bx^*)} \right| > \eta/2  \right) \nonumber\\
        &\leq 2 \exp \left\{ -\frac{\eta^2 \beta f(h^-(N)-1)}{8 (C_\mu^2 + C_{\text{sg}}^2)} \right\}.
    \end{align}
    
    RHS of \eqref{ineq:subgaussian_term2} converges to 0 as $N\to \infty$ because of $\lim_{N \to \infty} h^-(N) = \infty$ and $\lim_{c\to \infty} f(c) = \infty$. Then, the second term of \eqref{ineq:theorem_1} converges to $0$.
    
    \textbf{The third term.}

    Note that the randomness in the probability is reflected in $\cL_N$ and $\bar{Y}$. We have:
    \begin{align}
    \label{ineq:third_term}
        &\P\left(\bigcup_{L \in \cL_N} \{ \{ \dist(L, \bx^*) > \epsilon \} \cap \{ \bar{Y}_L \geq \mu(\bx^*) - \eta\}\}\right)\nonumber\\
        &= \E_{\cL_N} \left[ \P \left( \bigcup_{L \in \cL_N} \left\{ \{ \dist(L, \bx^*) > \epsilon \} \cap \{ \bar{Y}_L \geq \mu(\bx^*) - \eta \} \right\} | \cL_N \right)\right] \nonumber\\
        &\leq \E_{\cL_N} \left[ \sum_{L \in \cL_N}\P \left( \{ \dist(L, \bx^*) > \epsilon \} \cap \{ \bar{Y}_L \geq \mu(\bx^*) - \eta \} | L \right) \right],
    \end{align}
    where the equality follows from the law of total expectation, the inequality follows from Bonferroni's inequality.
    
    Taking $\eta = \delta(\epsilon)/2$, for any $\cL_N = \ell$, from the honesty that the estimation is independent of partitioning,
    \begin{align*}
        \sum_{L \in \ell}\P \left( \{ d(L, \bx^*) > \epsilon \} \cap \{ \bar{Y}_L \geq \mu(\bx^*) - \eta \} | L \right) &= \sum_{L \in \ell}\P \left( \{ d(L, \bx^*) > \epsilon \} \cap \{ \bar{Y}_L \geq \mu(\bx^*) - \eta \} \right)\\
        &\leq \sum_{L \in \ell} \P(\bar{Y}_L \geq \E[Y(L) ] + \delta(\epsilon)/2 )
    \end{align*}

     By honesty, $\E[\bar{Y}_{L}] = \E[Y(L)]$. From the Lemma \ref{lemma:subgaussian}, we know that for any $L$, $Y(L)$ is sub-Gaussian with variance proxy $C_\mu^2 + C_{\text{sg}}^2$. We have:
    \begin{align*}
        \sum_{L \in \ell} \P(\bar{Y}_L \geq \E[Y(L) ] + \delta(\epsilon)/2 ) &\leq \sum_{L \in \ell} \exp \left\{-  \frac{\beta f(h^-(N)-1) \delta(\epsilon)^2}{8(C_\mu^2 + C_{\text{sg}}^2)} \right\}\\
        &\leq 2^{h^+(N)} \exp \left\{-\frac{\beta \delta(\epsilon)^2}{8(C_\mu^2 + C_{\text{sg}}^2)} f(h^-(N)-1) \right\}
    \end{align*}
    where the first inequality follows from the Hoeffding inequality and the sample is greater than $\beta f(h^-(N)-1)$, the second inequality follows from the depth of leaf is upper bound by $h^+(N)$, which means that the leaf size is bounded above by $2^{h^+(N)}$.     
    
    Follow the assumption that $\lim_{N \to \infty}\frac{h^+(N)}{f(h^-(N))} = 0$, we know that the third term converges to 0.
\end{proof}

\subsection{Proof of Theorem \ref{theorem:optimal_gap_stochastic}}
\begin{proof}
    From the local smoothness of $\mu$, 
    \begin{align}
        \left| \bar{Y}(L_N^*) - \mu(\bx^*) \right| &\leq \left| \bar{Y}(L_N^*) - \E[Y( L_N (\bx^*) )|L_N^*] \right| + \left| \E[Y( L_N (\bx^*) )|L_N^*] - \mu(\bx^*) \right| \nonumber\\
        &\leq \left| \bar{Y}(L_N^*) - \E[Y( L_N (\bx^*) )|L_N^*] \right| + C \cdot \max_{\bx\in L_N^*} \| \bx - \bx^* \| \nonumber\\
        &\leq \left| \bar{Y}(L_N^*) - \E[Y( L_N (\bx^*) )|L_N^*] \right| + C \cdot (\diam(L_N^*) + \dist(L_N^*, \bx^*))\nonumber\\
        &= \left| \bar{Y}(L_N^*) - \E[Y( L_N (\bx^*) )|L_N^*] \right| + C \cdot \diam(L_N^*) + C \cdot \dist(L_N^*, \bx^*). \label{ineq:theorem_2}
    \end{align}

    By following the arguments presented in \eqref{ineq:subgaussian_term2}, the first term converges in probability to 0. Next, Lemma \ref{lemma:diam_depth} and $\lim_{N\to \infty} h^-(N) = \infty$ imply that the second term converges in probability to 0. Finally, according to Theorem \ref{theorem:sol_gap_stochastic_subgaussian}, the third term also converges in probability to 0. 

    Therefore, the overall error $\left| \bar{Y}(L_N^*) - \mu(\bx^*) \right|$ converges in probability to zero.
\end{proof}

\end{document}